\documentclass[MSNbibl,citesort,number,dvips]{arxbj}

\aid{0}
\volume{18}
\issue{1}
\pubyear{2012}
\firstpage{279}
\lastpage{289}
\doi{10.3150/10-BEJ337}

\makeatletter
\newtheorem{theorem}{Theorem}
\newtheorem{proposition}{Proposition}
\newtheorem{lemma}{Lemma}
\newtheorem{corollary}{Corollary}
\makeatother

\begin{document}
\begin{frontmatter}

\title{On the inclusion probabilities in some unequal probability sampling plans without replacement}
\runtitle{Sampling without replacement}

\begin{aug}
\author{\fnms{Yaming} \snm{Yu}\corref{}\ead[label=e1]{yamingy@uci.edu}}


\runauthor{Y. Yu}
\address{Department of Statistics,
University of California,
Irvine, CA 92697, USA.\\ \printead{e1}}
\end{aug}

\received{\smonth{6} \syear{2010}}

\begin{abstract}
Comparison results are obtained for the inclusion probabilities in some
unequal probability sampling plans without replacement.  For either
successive sampling or H\'{a}jek's rejective sampling, the larger the
sample size, the more uniform the inclusion probabilities in the sense
of majorization.  In particular, the inclusion probabilities are more
uniform than the drawing probabilities.  For the same sample size, and
given the same set of drawing probabilities, the inclusion
probabilities are more uniform for rejective sampling than for
successive sampling.  This last result confirms a conjecture of
H\'{a}jek (\textit{Sampling from a Finite Population} (1981) Dekker).  Results are
also presented in terms of the Kullback--Leibler divergence, showing
that the inclusion probabilities for successive sampling are more
proportional to the drawing probabilities.
\end{abstract}

\begin{keyword}
\kwd{conditional Poisson sampling}
\kwd{entropy}
\kwd{H\'{a}jek's conjecture}
\kwd{sampling
without replacement}
\kwd{stochastic orders}
\kwd{total positivity order}
\end{keyword}

\end{frontmatter}

\section{Introduction and main results}\label{sec1}

Consider a finite population indexed by $U=\{1, \ldots, N\}$.  Let
\mbox{$\alpha=(\alpha_1,\ldots, \alpha_N),\sum_{i=1}^N \alpha_i=1,$} denote a
set of drawing probabilities.  In H\'{a}jek's \cite{H64,H81}
\textit{rejective sampling}, independent draws are made with
probabilities according to the same $\alpha$ until a sample of size $n$
is obtained; whenever a duplicate appears, all draws are rejected and
the process restarts.  \textit{Successive sampling}, a closely related
scheme, makes the same independent draws except that whenever a
duplicate appears, only the current draw is rejected and needs to be
redrawn. Mathematically, rejective sampling is equivalent to
conditional Poisson sampling, that is, independent sampling for each
unit conditional on the sample size being $n$.  Conditional Poisson
sampling possesses a maximum entropy property, among other desirable
properties, and has received considerable attention; see Chen, Dempster
and Liu~\cite{CDL}, Berger \cite{YGB}, Traat, Bondesson and Meister
\cite{TBM}, Arratia, Goldstein and Langholz~\cite{AGL}, and
Qualit\'{e}\vadjust{\goodbreak}
\cite{Q}. It also has interesting applications to modeling how players
select lottery tickets \cite{SC}.  Successive sampling, on the other
hand, has connections to areas such as software reliability
\cite{Kauf}.

Unequal probability sampling may achieve considerable variance
reduction if the first-order inclusion probabilities are made
proportional to a suitable auxiliary variable.  For either rejective
sampling or successive sampling, however, the inclusion probabilities
are rather complicated and generally not proportional to the drawing
probabilities $\alpha$.  Thus relationships between the inclusion
probabilities and $\alpha$, either approximate or exact, are of
interest.  This work considers exact qualitative comparisons.  See
H\'{a}jek \cite{H64} and Ros\'{e}n~\cite{R1,R2,R3} for asymptotic
results.

Denote the inclusion probabilities for rejective sampling by
$\pi^{\mathrm{R}} =(\pi^{\mathrm{R}}_1, \ldots, \pi^{\mathrm{R}}_N)$
and those for successive sampling by $\pi^{\mathrm{S}}
=(\pi^{\mathrm{S}}_1, \ldots, \pi^{\mathrm{S}}_N)$. H\'{a}jek
\cite{H81}, page 97, conjectures the following inequalities based on
asymptotic considerations and numerical experience:
\[
\frac{\max \alpha_i}{\min \alpha_i} \geq \frac{\max
\pi^{\mathrm{S}}_i}{\min \pi^{\mathrm{S}}_i} \geq \frac{\max
\pi^{\mathrm{R}}_i}{\min \pi^{\mathrm{R}}_i}.
\]
Milbrodt \cite{Mil} proposes a strengthened conjecture,
\begin{eqnarray}
\label{mil1}
n \max \alpha_i & \geq& \max \pi^{\mathrm{S}}_i\geq \max \pi^{\mathrm{R}}_i,\\
\label{mil2} n \min \alpha_i & \leq& \min \pi^{\mathrm{S}}_i\leq \min
\pi^{\mathrm{R}}_i,
\end{eqnarray}
and partially resolves it by showing
\begin{eqnarray}
\label{wk1}
n \max \alpha_i &\geq &\max \pi^{\mathrm{S}}_i,\qquad  n\min\alpha_i\leq \min \pi^{\mathrm{S}}_i,\\
\nonumber n \max \alpha_i &\geq& \max \pi^{\mathrm{R}}_i,\qquad
n\min\alpha_i\leq \min \pi^{\mathrm{R}}_i.\nonumber
\end{eqnarray}
The inequalities (\ref{wk1}) are also obtained by Rao, Sengupta and
Sinha \cite{RSS}.  The inequalities $\max \pi^{\mathrm{S}}_i\geq \max
\pi^{\mathrm{R}}_i$ and $\min \pi^{\mathrm{S}}_i\leq \min
\pi^{\mathrm{R}}_i$ have remained open; see Milbrodt \cite{Mil} for
numerical illustrations. Roughly speaking, both H\'{a}jek's conjecture
and Milbrodt's strengthened version say that the drawing probabilities
are more variable than the inclusion probabilities for successive
sampling, which are themselves more variable than the inclusion
probabilities for rejective sampling.

Concerning successive sampling, Kochar and Korwar \cite{KK} obtain some
comparison results using the notion of majorization.  A real vector
$b=(b_1, \ldots, b_N)$ is said to majorize $a=(a_1, \ldots, a_N)$,
written as $a\prec b$, if
\begin{itemize}
\item $\sum_{i=1}^N a_i=\sum_{i=1}^N b_i$, and
\item $\sum_{i=k}^N
a_{(i)}\leq \sum_{i=k}^N b_{(i)}, k=2, \ldots, N,$ where $a_{(1)}\leq
\cdots\leq a_{(N)}$ and $b_{(1)}\leq \cdots\leq b_{(N)}$ are $(a_1,
\ldots, a_N)$ and $(b_1, \ldots, b_N)$ arranged in increasing order,
respectively.
\end{itemize}
Kochar and Korwar \cite{KK} show that
\begin{equation}
\label{kk} n^{-1} \pi^{\mathrm{S}} \prec \alpha,
\end{equation}
which strengthens (\ref{wk1}).  In general, majorization is a strong
form of variability ordering.  For example, $a\prec b$ implies that
$\sum_i \phi(a_i) \leq \sum_i \phi(b_i)$ for any convex function
$\phi$.  See Marshall and Olkin \cite{MO} for further properties and
various applications of majorization.\vadjust{\goodbreak}

This note presents some majorization results that refine previous work.
As a consequence, we prove Milbrodt's strengthening of H\'{a}jek's
conjecture.  Our main results are summarized as follows.

\begin{theorem}
\label{thm1} Given the drawing probabilities $\alpha$, let
$\pi^{\mathrm{R}}(n)$ (resp., $\pi^{\mathrm{S}}(n)$) denote the
first-order inclusion probabilities for rejective sampling (resp.,
successive sampling) with sample size \mbox{$n\leq N$}.  Define the
``inclusion probabilities per draw'' as $p^{\mathrm{R}}(n)\equiv n^{-1}
\pi^{\mathrm{R}}(n)$ and $p^{\mathrm{S}}(n)\equiv n^{-1}
\pi^{\mathrm{S}}(n)$. Then we have
\begin{eqnarray}
\label{majR}
(N^{-1},\ldots, N^{-1})&\equiv& p^{\mathrm{R}}(N) \prec \cdots \prec  p^{\mathrm{R}}(n)\prec \cdots \prec p^{\mathrm{R}}(1)\equiv \alpha,\\
\label{majS} (N^{-1},\ldots, N^{-1})&\equiv& p^{\mathrm{S}}(N) \prec
\cdots \prec  p^{\mathrm{S}}(n)\prec \cdots \prec
p^{\mathrm{S}}(1)\equiv \alpha.
\end{eqnarray}
Moreover,
\begin{equation}
\label{piRS} \pi^{\mathrm{R}}(n) \prec \pi^{\mathrm{S}}(n).
\end{equation}
\end{theorem}

The ordering chains (\ref{majR}) and (\ref{majS}) are intuitively
appealing.  Given a set of drawing probabilities, larger sample sizes
lead to inclusion probabilities that are more uniform for either
rejective sampling or successive sampling.  Moreover, (\ref{piRS}) says
that with the same sample size, the inclusion probabilities are more
uniform for rejective sampling than for successive sampling.  It is
easy to see that (\ref{majR})--(\ref{piRS}) together imply Milbrodt's
\cite{Mil} conjecture, that is, (\ref{mil1}) and (\ref{mil2}).

We prove (\ref{majR}) and (\ref{piRS}) in Section~\ref{sec2} using a combination
of analytic and probabilistic techniques.  A key tool in resolving
(\ref{piRS}) is the likelihood ratio order between multivariate
densities \cite{KR80}.  A proof of (\ref{majS}), which slightly extends
that of (\ref{kk}), is included for completeness.\looseness=-1

The\vspace*{1pt} Shannon entropy is sometimes used to measure how uniform a
distribution is.  It is defined as $H(p)=-\sum_{i=1}^N p_i\log p_i$ for
a probability vector $p=(p_1,\ldots, p_N)$.  By convention \mbox{$0\log 0=0$}.
It is well known that $p\prec q$ implies $H(q)\leq H(p)$.  See Cover
and Thomas \cite{Cover}, Chapter~2, for further properties of this
fundamental quantity.  We note the following direct consequence of
Theorem~\ref{thm1}.

\begin{corollary}
In the setting of Theorem \ref{thm1},
\begin{eqnarray*}
\log N&\equiv& H(p^{\mathrm{R}}(N)) \geq  \cdots \geq  H(p^{\mathrm{R}}(n))\geq \cdots \geq H(p^{\mathrm{R}}(1))\equiv H(\alpha),\\
\log N&\equiv& H(p^{\mathrm{S}}(N)) \geq  \cdots \geq  H(p^{\mathrm{S}}(n))\geq \cdots \geq H(p^{\mathrm{S}}(1))\equiv H(\alpha),\\
H(p^{\mathrm{R}}(n)) &\geq&  H(p^{\mathrm{S}}(n)).
\end{eqnarray*}
\end{corollary}

Inequalities are also obtained in terms of the Kullback--Leibler
divergence, which is defined as
\[
D(p\|q)= \sum_{i=1}^N p_i \log
\frac{p_i}{q_i}
\]
for two probability vectors $p=(p_1,\ldots, p_N)$ and $q=(q_1,\ldots,
q_N)$.  By convention $x\log(x/0)=\infty$ for $x>0$ and $0\log (0/x)=0$
for $x\geq 0$.  A basic property is $D(p\|q)> 0$ unless $p=q$.  We
shall use $D(p\|q)$ purely as a discrepancy measure between probability
vectors without referring to its information-theoretic significance.

\begin{theorem}
\label{thm2} In the setting of Theorem \ref{thm1}, let $1\leq l< m<
n\leq N$.  Then we have
\begin{eqnarray}
\label{r1}
D(p^{\mathrm{R}}(l)\|p^{\mathrm{R}}(n)) &\geq& D(p^{\mathrm{R}}(l)\|p^{\mathrm{R}}(m)) + D(p^{\mathrm{R}}(m)\| p^{\mathrm{R}}(n)),\\
\label{r2}
D(p^{\mathrm{R}}(n)\|p^{\mathrm{R}}(l)) &\geq& D(p^{\mathrm{R}}(m)\|p^{\mathrm{R}}(l)) + D(p^{\mathrm{R}}(n)\| p^{\mathrm{R}}(m)),\\
\label{ss}
D(p^{\mathrm{S}}(n)\|\alpha) &\geq& D(p^{\mathrm{S}}(n)\|p^{\mathrm{S}}(m)) + D(p^{\mathrm{S}}(m)\| \alpha),\\
\label{rs} D(p^{\mathrm{R}}(n)\|\alpha) &\geq&
D(p^{\mathrm{R}}(n)\|p^{\mathrm{S}}(n)) + D(p^{\mathrm{S}}(n)\|
\alpha).
\end{eqnarray}
\end{theorem}

A number of results can be deduced from these (reverse) triangle
inequalities.  For example, from (\ref{r1}) and (\ref{r2}) we obtain
$D(p^{\mathrm{R}}(m+1)\|\alpha) \geq D(p^{\mathrm{R}}(m)\|\alpha)$ and
$D(\alpha\|p^{\mathrm{R}}(m+1)) \geq D(\alpha\|p^{\mathrm{R}}(m))$,
showing that, for rejective sampling, the larger the sample size, the
more distorted the inclusion probabilities become as compared with the
drawing probabilities.  Similarly, from (\ref{ss}) we obtain
$D(p^{\mathrm{S}}(m+1)\|\alpha) \geq D(p^{\mathrm{S}}(m)\|\alpha)$.
From (\ref{rs}) we obtain
\begin{equation}
\label{bd} D(p^{\mathrm{R}}(n)\|\alpha) \geq
D(p^{\mathrm{S}}(n)\|\alpha).
\end{equation}
That is, for fixed $n$, the inclusion probabilities for successive
sampling (rather than for rejective sampling) are more proportional to
the drawing probabilities.  The inequality~(\ref{bd}) may be used to
compute an upper bound on $D(p^{\mathrm{S}}(n)\|\alpha)$ because, while
$p^{\mathrm{R}}(n)$ can be calculated from $\alpha$ efficiently using a
recursive formula (see \cite{CDL}), numerical calculation of~$p^{\mathrm{S}}(n)$ is considerably more difficult.

The inequalities in Theorem~\ref{thm2} resemble the reverse triangle
inequalities of Yu \cite{Yu}.  Our results here concern the
majorization ordering and may be regarded as first-order results; those
in Yu \cite{Yu} use relative log-concavity and are second order.  For
related entropy and divergence comparison results, see Karlin and
Rinott \cite{KR81}, Johnson \cite{J} and Yu \cite{Yu2,Yu3,Yu4}.

The proof of Theorem~\ref{thm2} builds on Theorem~\ref{thm1} and is
presented in Section~\ref{sec3}.

\section{\texorpdfstring{Proof of Theorem~\protect\ref{thm1}}{Proof of Theorem 1}}\label{sec2}
Let $e_k(\cdot)$ denote the $k$th elementary symmetric function, that
is,
\[
e_k(\beta) = \sum_{1\leq j_1<\cdots < j_k\leq m} \beta_{j_1}
\cdots \beta_{j_k},\qquad \beta\equiv (\beta_1, \ldots, \beta_m).
\]
 By
convention, $e_0(\beta)\equiv 1$ and $e_k(\beta)=0$ if $k<0$ or $k>m$.
For a rejective sample of size~$n$, the probability that unit $i$ is
included can be expressed as
\begin{equation}
\label{incR} \pi^{\mathrm{R}}_i(n) = \frac{\alpha_i
e_{n-1}(\alpha_{-i})}{ e_n(\alpha)},
\end{equation}
where
\[
\alpha_{-i}\equiv (\alpha_1, \ldots, \alpha_{i-1}, \alpha_{i+1}, \ldots, \alpha_N).
\]
The notation $\alpha_{-i, -j}$ (leave-two-out) is defined similarly. It
is immediate that \mbox{$\alpha_i \leq \alpha_j$} implies
$\pi^{\mathrm{R}}_i(n) \leq \pi^{\mathrm{R}}_j(n)$.  Henceforth, we
assume $\alpha_1\geq \cdots \geq \alpha_N>0$ without loss of
generality.\looseness=1

The following Lemma~\ref{lem1} is needed in the proof of~(\ref{majR}).

\begin{lemma}
\label{lem1}
Suppose probability vectors $p=(p_1, \ldots, p_N)$ and $q
=(q_1,\ldots, q_N)$ satisfy
\[
p_1\geq \cdots\geq p_N>0,\qquad \frac{q_1}{p_1}\geq \cdots \geq \frac{q_N}{p_N}.
\]
Then $p\prec q.$
\end{lemma}

\begin{pf}
For $1\leq k< N$ we have
\[
\frac{\sum_{i=1}^k q_i}{\sum_{i=1}^k p_i} \geq \frac{q_k}{p_k}\geq \frac{q_{k+1}}{p_{k+1}},
\]
which yields
\[
\frac{\sum_{i=1}^k q_i}{\sum_{i=1}^k p_i} \geq
\frac{\sum_{i=1}^{k+1} q_i}{\sum_{i=1}^{k+1} p_i}\geq \cdots \geq \frac{\sum_{i=1}^N q_i}{\sum_{i=1}^N p_i} =1.
\]
Hence $p\prec q$ by definition (the conditions imply $q_1\geq \cdots
\geq q_N$).
\end{pf}

\begin{pf*}{Proof of (\ref{majR})}
Let $p\equiv p^{\mathrm{R}}(n+1)$ and $q\equiv p^{\mathrm{R}}(n)$. Note
that $\sum_{i=1}^N p_i =\sum_{i=1}^N q_i =1$.  Since $\alpha_1\geq
\cdots \geq \alpha_N$, we have $p_1\geq\cdots\geq p_N$.  The desired
relation $p\prec q$ would follow from Lemma~\ref{lem1}, if we can show
that $q_1/p_1\geq \cdots \geq q_N/p_N$, or, equivalently,
$\pi^{\mathrm{R}}_k(n)/\pi^{\mathrm{R}}_k(n+1)\geq
\pi^{\mathrm{R}}_{k+1}(n)/\pi^{\mathrm{R}}_{k+1}(n+1)$ for $1\leq k
<N$. The case $N=2$ is trivial.  Otherwise we have
\[
\pi^{\mathrm{R}}_k(n) = \frac{\alpha_k e_{n-1}(\alpha_{-k})}{ e_n(\alpha)}
=\alpha_k \frac{\alpha_{k+1} e_{n-2} (\tilde{\alpha}) + e_{n-1}(\tilde{\alpha})}{e_n(\alpha)}, \qquad\tilde{\alpha}\equiv \alpha_{-k, -(k+1)}.
\]
Thus
\begin{equation}
\label{rk} \frac{\pi^{\mathrm{R}}_k(n)}{\pi^{\mathrm{R}}_k(n+1)} =
\frac{e_{n+1}(\alpha)}{e_n(\alpha)} f(\alpha_{k+1}),
\end{equation}
where
\[
f(x) =  \frac{x e_{n-2} (\tilde{\alpha}) + e_{n-1}(\tilde{\alpha})}{x e_{n-1} (\tilde{\alpha}) + e_n(\tilde{\alpha})}.
\]
Similarly
\begin{equation}
\label{rk+}
\frac{\pi^{\mathrm{R}}_{k+1}(n)}{\pi^{\mathrm{R}}_{k+1}(n+1)} =
\frac{e_{n+1}(\alpha)}{e_n(\alpha)} f(\alpha_k).
\end{equation}
We have
\[
f'(x) = \frac{e_{n-2} (\tilde{\alpha}) e_n(\tilde{\alpha}) - e_{n-1}^2(\tilde{\alpha})}{[x e_{n-1} (\tilde{\alpha}) + e_n(\tilde{\alpha})]^2} < 0,
\]
where the inequality follows from Newton's inequalities \cite{HLP}, page 52.  That is, $f(x)$ decreases in~$x$.  Because $\alpha_{k+1} \leq \alpha_{k}$, we deduce the inequality
\[
\frac{\pi^{\mathrm{R}}_k(n)}{\pi^{\mathrm{R}}_k(n+1)}\geq \frac{\pi^{\mathrm{R}}_{k+1}(n)}{\pi^{\mathrm{R}}_{k+1}(n+1)}
\]
from (\ref{rk}) and (\ref{rk+}).
\end{pf*}

The proof of (\ref{majS}) slightly extends and simplifies the arguments
of Kochar and Korwar~\cite{KK}.

\begin{pf*}{Proof of (\ref{majS})}
Let $S_1, S_2, \ldots\in \{1,\ldots,N\}$ be a sequence of draws
retained in successive sampling.  It is well known that the inclusion
probabilities and the drawing probabilities are ordered in the same
way, that is,
\begin{equation}
\label{sorder} p^{\mathrm{S}}_1(n)\geq \cdots\geq p^{\mathrm{S}}_N(n),\qquad
1\leq n\leq N
\end{equation}
(see \cite{Mil}).  For $1\leq k\leq N$ we have
\begin{eqnarray*}
&&\Pr(S_n \leq k)-\Pr(S_{n+1}\leq k) \\
&&\quad= \Pr(S_n \leq k, S_{n+1}>k) -\Pr(S_n>k, S_{n+1}\leq k)\\
&&\quad= \sum_{k_1\leq k, k_2>k} \mathbf{E} [\Pr(S_n =k_1, S_{n+1}=k_2 |S_1,\ldots, S_{n-1})\\
&&\qquad\hphantom{\sum_{k_1\leq k, k_2>k} \mathbf{E} [}{}   -\Pr(S_n=k_2, S_{n+1}=k_1 |S_1,\ldots, S_{n-1})],
\end{eqnarray*}
where the expectation is with respect to $S_1, \ldots, S_{n-1}$.
Because $\alpha_i$ decreases in $i$, it is easy to show that $k_1<k_2$
implies
\[
\Pr(S_n =k_1, S_{n+1}=k_2 |S_1,\ldots, S_{n-1}) \geq \Pr(S_n=k_2, S_{n+1}=k_1|S_1,\ldots, S_{n-1}).
\]
Hence $\Pr(S_n \leq k)\geq \Pr(S_{n+1}\leq k)$ for all $1\leq n< N$.
This is proved by Kochar and Korwar \cite{KK} (see their Lemma 3.2) using
a slightly more complicated argument.  It follows that
\begin{eqnarray*}
\sum_{i=1}^k p^{\mathrm{S}}_i(n) &=& n^{-1} \sum_{j=1}^n \Pr(S_j \leq k)\\
&\geq& (n+1)^{-1} \sum_{j=1}^{n+1} \Pr(S_j\leq k)\\
&=& \sum_{i=1}^k p^{\mathrm{S}}_i(n+1),
\end{eqnarray*}
which proves (\ref{majS}) in view of (\ref{sorder}).\vadjust{\goodbreak}
\end{pf*}

To prove (\ref{piRS}), we recall the multivariate likelihood ratio
order, also known as the total positivity order (Karlin and Rinott
\cite{KR80}, Rinott and Scarsini \cite{RS}, Shaked and Shanthikumar \cite{SS07}, Chapter~6).
Consider the product space $\mathcal{X}=\{1, \ldots, N\}^n$.  For
$x=(x_1, \ldots, x_n)\in \mathcal{X}$ and $y=(y_1, \ldots, y_n)\in
\mathcal{X}$, write
\[
x\vee y = (\max\{x_1, y_1\}, \ldots, \max\{x_n, y_n\}), \qquad x\wedge y =(\min\{x_1, y_1\}, \ldots, \min\{x_n, y_n\}).
\]
Let $f$ and $g$ be density functions on $\mathcal{X}$.  Then $f$ is
said to be no smaller than $g$ in the (multivariate) likelihood ratio
order, written as $f \geq_{\mathrm{lr}} g$, if
\[
f(x)g(y) \leq f(x\vee y) g(x\wedge y),\qquad x, y\in \mathcal{X}.
\]
This generalizes the univariate likelihood ratio order, which requires
that the ratio of two univariate densities is a monotone function.

A useful property of the likelihood ratio order is that it implies the
usual stochastic order.  That is, if $X$ and $Y$ are random vectors
taking values in $\mathcal{X}$, and $X\geq_{\mathrm{lr}} Y$ (we use the
notation $\geq_{\mathrm{lr}}$ with the random variables as well as
their densities), then $\mathbf{E}\phi(X)\geq \mathbf{E}\phi(Y)$ for
any coordinatewise increasing function $\phi$.  In particular, each
coordinate of $X$ is no smaller than the corresponding coordinate of
$Y$ in the usual stochastic order.  Further properties of
$\geq_{\mathrm{lr}}$ include closure under marginalization; see Karlin
and Rinott \cite{KR80} and Shaked and Shanthikumar \cite{SS07}, Chapter~6.

\begin{pf*}{Proof of (\ref{piRS})}
Recall that $\pi^{\mathrm{R}}_1(n)\geq \cdots \geq
\pi^{\mathrm{R}}_N(n).$ By definition, (\ref{piRS}) is proved if we can
show
\begin{equation}
\label{goal} \sum_{i=1}^k \pi^{\mathrm{S}}_i(n) \geq \sum_{i=1}^k
\pi^{\mathrm{R}}_i(n),\qquad k=1,\ldots, N-1.
\end{equation}

Let $X\equiv (X_1, \ldots, X_n)$ (resp., $Y\equiv (Y_1, \ldots,
Y_n)$) denote the unit indices arranged in increasing order of a sample
of size $n$ obtained by rejective sampling (resp., successive
sampling).  That is, $X$ and $Y$ take values in \mbox{$\Omega \equiv
\{(x_1,\ldots, x_n)\in \mathcal{X}\dvt 1\leq x_1<\cdots<x_n\leq N\}.$} Then
an unnormalized density of $X$ is
\[
f(x) = \alpha_{x_1} \ldots \alpha_{x_n},\qquad x=(x_1,\ldots, x_n)\in \Omega,
\]
and the density of $Y$ can be written as
\begin{eqnarray*}
g(y) &=& \sum_{\sigma\in \operatorname{Perm}(y)}
\alpha_{\sigma_1}\frac{\alpha_{\sigma_2}}{1-\alpha_{\sigma_1}}
\frac{\alpha_{\sigma_3}}{1-\alpha_{\sigma_1}-\alpha_{\sigma_2}}\cdots \frac{\alpha_{\sigma_n}}{1-\sum_{j=1}^{n-1} \alpha_{\sigma_j}}\\[3pt]
&=&  \sum_{\sigma\in \operatorname{Perm}(y)} \frac{\alpha_{y_1}\cdots
\alpha_{y_n}}{(1-\alpha_{\sigma_1})
(1-\alpha_{\sigma_1}-\alpha_{\sigma_2})\cdots (1-\sum_{j=1}^{n-1}
\alpha_{\sigma_j})},\qquad y=(y_1, \ldots, y_n)\in \Omega,
\end{eqnarray*}
where $\sigma=(\sigma_1, \ldots, \sigma_n)$ and $\operatorname{Perm}(y)$ denotes
the set of vectors obtained by permuting the coordinates of $y$.  Note
that, for $x, y\in \Omega$ we have $x\vee y\in \Omega$ and $x\wedge y
\in \Omega$.  Moreover, for $x, y\in \Omega$,
\begin{eqnarray*}
f(x)g(y) &= &\sum_{\sigma\in \operatorname{Perm}(y)} \frac{\alpha_{x_1}\cdots
\alpha_{x_n} \alpha_{y_1}\cdots \alpha_{y_n}}{(1-\alpha_{\sigma_1})
(1-\alpha_{\sigma_1}-\alpha_{\sigma_2})\cdots (1-\sum_{j=1}^{n-1} \alpha_{\sigma_j})}\\[3pt]
& \leq& \sum_{\sigma\in \operatorname{Perm}(x\wedge y)} \frac{\alpha_{x_1}\cdots
\alpha_{x_n} \alpha_{y_1}\cdots \alpha_{y_n}}{(1-\alpha_{\sigma_1})
(1-\alpha_{\sigma_1}-\alpha_{\sigma_2})\cdots (1-\sum_{j=1}^{n-1} \alpha_{\sigma_j})}\\[3pt]
& =& f(x\vee y) g(x\wedge y),
\end{eqnarray*}
where the inequality holds because $\alpha_i$ decreases in $i$ and,
under an obvious bijection, each element in $\operatorname{Perm}(y)$ is at least
as large as its counterpart in $\operatorname{Perm}(x\wedge y)$.  Thus
$X\geq_{\mathrm{lr}} Y$.  It follows that
\[
\Pr(X_j\leq k)\leq \Pr(Y_j\leq k),\qquad j=1,\ldots, n, k=1,\ldots, N.
\]
That is, $X_j$ is no smaller than $Y_j$ in the usual stochastic order.
We have, for $1\leq k\leq N$,
\begin{eqnarray*}
\sum_{i=1}^k \pi^{\mathrm{R}}_i(n) &=&\sum_{i=1}^k \sum_{j=1}^n \Pr(X_j=i)\\
& =& \sum_{j=1}^n \Pr(X_j\leq k) \\
&\leq &\sum_{j=1}^n \Pr(Y_j\leq k)\\
&=&\sum_{i=1}^k \sum_{j=1}^n \Pr(Y_j=i)\\
& =&\sum_{i=1}^k \pi^{\mathrm{S}}_i(n).
\end{eqnarray*}
Thus (\ref{goal}) holds, and the proof is complete.
\end{pf*}

\section{\texorpdfstring{Proof of Theorem~\protect\ref{thm2}}{Proof of Theorem 2}}\label{sec3}

The following Lemma~\ref{lem2} is key to the proof of
Theorem~\ref{thm2}.

\begin{lemma}
\label{lem2} Let $p=(p_1,\ldots, p_N), q=(q_1,\ldots, q_N)$ and
$r=(r_1,\ldots, r_N)$ be probability vectors with all positive
coordinates.  If either \textup{(a)} $q\prec p, p_1\geq \cdots \geq p_N$, and
$q_1/r_1\geq \cdots\geq q_N/r_N,$ or \textup{(b)} $p\prec q, q_1\geq \cdots \geq
q_N$, and $q_1/r_1\leq \cdots\leq q_N/r_N,$ then
\[
D(p\|r)\geq D(p\|q)+D(q\|r).
\]
\end{lemma}

\begin{pf}
Let us assume (a).  Case (b) is similar.  We have
\begin{eqnarray}\label{abel}
D(p\|r)-D(p\|q)-D(q\|r) &= &\sum_{i=1}^N (p_i-q_i) \log\frac{q_i}{r_i}\nonumber\\[3pt]
& =& \sum_{i=1}^{N-1} \Biggl(\sum_{j=1}^i p_j -\sum_{j=1}^i q_j\Biggr) \biggl(\log\frac{q_i}{r_i}-\log\frac{q_{i+1}}{r_{i+1}}\biggr)\\[3pt]
 & \geq& 0,\nonumber
\end{eqnarray}
where the first equality follows from the definition of the
Kullback--Leibler divergence, the second equality holds by summation by
parts, and the inequality holds because $q_i/r_i$ decreases in $i$ and
$q\prec p$, and hence both parentheses in (\ref{abel}) are
non-negative.
\end{pf}

As in Section~\ref{sec2}, in the proofs of (\ref{r1})--(\ref{rs}) we assume
$\alpha_1\geq \cdots \geq \alpha_N>0$.

\begin{pf*}{Proof of (\ref{r1}) and (\ref{r2})}
Let $p\equiv p^{\mathrm{R}}(l), q\equiv p^{\mathrm{R}}(m), r\equiv
p^{\mathrm{R}}(n)$.  Then $p_1\geq \cdots \geq p_N$.  Since $l<m$ we
have $q\prec p$ by (\ref{majR}).  From the proof of (\ref{majR}) we
know that $q_1/r_1\geq \cdots \geq q_N/r_N.$  Thus (\ref{r1}) follows
from Lemma~\ref{lem2}, Case (a).  The proof of (\ref{r2}) is similar.
\end{pf*}

To prove (\ref{ss}) and (\ref{rs}) we need the following result.

\begin{proposition}
\label{prop1}
The ratio $p^{\mathrm{S}}_i(n)/\alpha_i, i=1,\ldots, N,$
increases in $i$ for each $n\leq N$.
\end{proposition}

\begin{pf}
Let $\pi_{i, k}$ denote the probability that the $k$th distinct draw in
successive sampling is unit $i$.  Then
$p^{\mathrm{S}}_i(n)=n^{-1}\sum_{k=0}^{n-1} \pi_{i,k+1}$.  It suffices
to show that $\pi_{i,k+1}/\alpha_i$ increases in $i$ for each $k$.  Let
us assume $k\geq 1$ and define the index set
\[
\Omega(i)=\{(j_1, \ldots, j_{k})\dvt 1\leq j_l\leq N, j_l\neq i, 1\leq l\leq k, \mbox{ and $j_l$ are distinct}\}.
\]
Then we have
\begin{equation}
\label{longsum} \frac{\pi_{i, k+1}}{\alpha_i}=\sum_{(j_1,\ldots,
j_{k})\in \Omega(i)} \alpha_{j_1} \frac{\alpha_{j_2}}{1-\alpha_{j_1}}
\cdots \frac{\alpha_{j_{k}}}{1-\sum_{l=1}^{k-1} \alpha_{j_l}}
\biggl(\frac{1}{1-\sum_{l=1}^{k} \alpha_{j_l}}\biggr).
\end{equation}
The summand is a decreasing function in $(j_1,\ldots,j_k)$, since
$\alpha_j$ decreases in $j$.  Consider a~mapping $\Omega(i) \to
\Omega(i+1)$ that sends $(j_1, \ldots, j_{k})\in \Omega(i)$ to $(j^*_1,
\ldots, j^*_k)\in \Omega(i+1)$ as follows.  For $l=1,\ldots, k$, if
$j_l\neq i+1,$  let $j^*_l=j_l$; otherwise let $j^*_l=i$.  It is easy
to see that this mapping is well defined and is a bijection.  Note that
$j^*_l\leq j_l$.  Hence the right-hand side of (\ref{longsum})
increases if we replace the summation index $\Omega(i)$ by
$\Omega(i+1)$.  That is, $\pi_{i, k+1}/\alpha_i$ increases in $i$, as
required.
\end{pf}

\begin{pf*}{Proof of (\ref{ss}) and (\ref{rs})}
Let $p\equiv p^{\mathrm{R}}(n), q\equiv p^{\mathrm{S}}(n)$ and $r\equiv
\alpha$. By Proposition \ref{prop1}, $q_1/r_1\leq \cdots \leq q_N/r_N$.
By (\ref{piRS}) we have $p\prec q$.  Thus (\ref{rs}) follows from
Lemma~\ref{lem2}, Case (b).  The proof of (\ref{ss}) is similar.
\end{pf*}


\printhistory

\end{document}